\documentclass[11pt,reqno]{amsart}
\usepackage{amsmath,amsthm,amsfonts,amssymb,amscd,latexsym}
\usepackage[all]{xy}

\def\F{\mathbb{F}}

\def\su{{\subseteq}}

\def\s{{\bf s}}

\def\charac{{\rm char}}

\def\ud{{\tilde{\delta}}}

\def\Proof{\noindent{\sl Proof.}\ }
\def\qed{{\hfill $\Box$ \medbreak}}

\newtheorem{defi}{Definition}[section]
\newtheorem{thm}[defi]{Theorem}
\newtheorem{lem}[defi]{Lemma}
\newtheorem{cor}[defi]{Corollary}

\newtheorem{prop}[defi]{Proposition}
\newtheorem{rem}[defi]{Remark}

\begin{document}

\title[Solvability of Poisson algebras]{Solvability of Poisson algebras}
\author{\textsc{Salvatore Siciliano}}
\address{Dipartimento di Matematica e Fisica ``Ennio De Giorgi", Universit\`{a} del Salento,
Via Provinciale Lecce--Arnesano, 73100--Lecce, Italy}
\email{salvatore.siciliano@unisalento.it}

\author{\textsc{Hamid Usefi}}
\address{Department of Mathematics and Statistics,
Memorial University of Newfoundland,
St. John's, NL,
Canada, 
A1C 5S7}
\email{usefi@mun.ca}

\begin{abstract}  Let $P$ be a Poisson algebra with a Lie bracket $\{, \}$ over a field $\F$ of characteristic $p\geq 0$. In this paper, the Lie structure of $P$ is investigated. In particular, if $P$ is solvable with respect to its Lie bracket, then we prove that  the Poisson ideal $\mathcal{J}$  of $P$ generated by all elements $\{\{\{x_1, x_2\}, \{x_3, x_4\}\}, x_5\}$ with $x_1,\ldots ,x_5 \in P$ is associative nilpotent of index bounded by a function of the derived length of $P$. We use this result to further prove that if $P$ is solvable  and $p\neq 2$, then the Poisson ideal $\{P,P\}P$ is nil.  
\end{abstract}
\date{\today}

  \maketitle

\section{Introduction} 

Since the theory of polynomial identities (PI) has proven to be  extremely useful in the variety of associative, Lie, or Jordan algebras, in the past decade there has been great interest to develop analogues theories for Poisson algebras.
Recall that a Poisson algebra $P$ is a commutative associative algebra with a unity equipped with a Lie bracket $\{, \}$ that satisfies the Leibniz rule:
\begin{align*}
\{a \cdot b, c\} = a \cdot \{b, c\} + b \cdot \{a, c\}, \quad a, b, c \in P.
\end{align*}
Poisson algebras  naturally arise in different areas of algebra, topology and mathematical physics and have received a considerable attention over the years.
One of the main initiatives in the domain of PI Poisson algebras was taken by Farkas \cite{F1, F2} who introduced the analogues of standard polynomial identities. This work was further picked up in \cite{MPR} to study codimension growth in characteristic zero and prove that the tensor product of PI Poisson algebras is again a PI-algebra. Moreover, in \cite{GP} Giambruno and Petrogradsky established when the symmetric Poisson algebra $S(L)$ or the  truncated symmetric Poisson algebra $\s(L)$ of a restricted Lie algebra $L$ satisfies a nontrivial multilinear Poisson identical relation. More recently,  in \cite{MP}, Monteiro Alves and Petrogradsky focused on Lie identities of $S(L)$ and $\s(L)$, establishing in particular when these Poisson algebras are Lie nilpotent or solvable in odd characteristic. Further developments of these topics have been recently carried out by the first author in \cite{S}.

In this paper, we investigate the Lie identities of an arbitrary Poisson algebra $P$ and study the relationship between the Lie structure  and the associative structure of $P$. Motivation for these problems also arises from analogous problems for rings.
An important result for  Lie solvable varieties of associative algebras was proved independently by 
Sharma and Srivastava \cite{SS} and Smirnov and Zalesskii \cite{SZ}.
Let $R$ be an associative algebra and consider the Lie product on $R$ given by $[x,y]=xy-yx$ for all $x,y\in R$. Denote by $\mathcal{J}$ the ideal of $R$ generated by all elements 
$[[[x_1, x_2], [x_3, x_4]], x_5]$ with $x_1,\ldots ,x_5 \in R$. It is proved in \cite{SS} and \cite{SZ} that if $R$ is Lie solvable, then $\mathcal{J}$ is associative nilpotent. This result has been used numerously and  proven to be very useful (see e.g. \cite{AS, R, SU1, SU2}). The primary goal of this paper is to prove a similar statement for Poisson algebras. In Theorem \ref{S-Z}, we prove that if $P$ is a solvable Poisson algebra of derived length $n$, then the Poisson ideal generated by all elements $\{\{\{x_1,x_2\},\{x_3,x_4\}\},x_5\}$ is associative nilpotent of index bounded by a function of $n$. We also mention that an application of this result in a forthcoming paper allows to settle the solvability problem of $S(L)$ and $\s(L)$ in characteristic 2 posed in \cite[\S 5.3]{MP}. 

In our second main result, we investigate the extent to which the ideal  $\{P, P\}P$ is  nil  (of bounded index) in case $P$ is solvable. 
Jennings \cite{J} proved that if $R$ is a finitely generated Lie nilpotent ring, then the ideal $[R, R]R$  is nilpotent. In \cite{R}, for an associative algebra $A$ over a field of characteristic $p>0$, Riley showed that $[A, A]A$ is nil of bounded index whenever $A$ is Lie nilpotent, or Lie solvable and $p>2$. Now, in characteristic zero, it follows from \cite[Theorem 7.2]{MPR} that if a Poisson algebra $P$ is solvable, then $\{P,P\}P$ is nil. In Theorem \ref{solvablenil}, we strengthen this result by proving that if $P$ is solvable over a field of characteristic $p\neq 2$, then the Poisson ideal $\{P,P\}P$ is nil. Furthermore, if  $P$ is solvable and  $p\geq  3$ or $P$ is Lie nilpotent and  $p\geq  2$, then $\{P,P\}P$ is nil of bounded index. This generalizes the main result in \cite{MP} about the Lie structure of $S(L)$   to all reduced Poisson algebras.

\section{Definitions and notation}
  Let $P$ be a Poisson algebra over a field $\F$. 
 The Poisson brackets are left-normed: $\{x_1,\ldots,x_n\}=\{\{x_1,\ldots,x_{n-1}\},x_n\}$, $n\geq 1$.
We will denote by $\langle S \rangle_\F$ the subspace spanned by a subset $S$ of $P$. We  use the symbol  $Z_P(P)$ for the Poisson center of $P$ (that is, the center of $P$ as a Lie algebra). When one deals with ideals of Poisson algebras or their properties, it is important
to distinguish which operation is considered. Our convention is
that the term \emph{ideal} refers to the associative multiplication, \emph{Lie ideal} refers to the Lie bracket, and \emph{Poisson ideal} refers to
both. 

The terms of the derived series of a Lie ideal $L$ of $P$ are defined by $\delta_0(L) = L$ and $\delta_{n}(L) = \{\delta_{n-1}(L), \delta_{n-1}(L)\}$ for $n>0$. Moreover, we denote by $\gamma_n(L)$ ($n\geq 1$) the terms of the descending central series of $L$. In particular, $\gamma_2(L)=\delta_1(L)$  is the derived subalgebra of $L$. We say that $P$ is \emph{ Lie nilpotent} if $P$ is nilpotent as a Lie algebra. The upper Lie power series of $P$ is the chain of Poisson ideals of $P$ defined by $P^{(1)}=P$ and $P^{(n)}=\{P^{(n-1)},P^{(n-1)}\}P$ for every $n>1$. The  Poisson algebra $P$ is said to be \emph{strongly Lie nilpotent} of class $c$ if $P^{(c+1)}=0$ and $P^{(c)}\neq 0$.  
  One says that the Poisson algebra $P$ is \emph{solvable} if $P$ is solvable as a Lie algebra. In this case, the minimal $n$ such that   $\delta_{n}(P)=0$ is called the \emph{derived length} of $P$. 
	The upper derived series of $P$ is defined by setting $\ud_0(P)=P$ and $\ud_n(P)=\{\ud_{n-1}(P),\ud_{n-1}(P)\} P$ for every $n>0$. Note that $\ud_n(P)$ is a Poisson ideal of $P$ for every $n$. The Poisson algebra $P$ is said to be \emph{strongly solvable} if $\ud_n(P)=0$, for some $n$. In this case, the minimal $n$ with such a property is called the \emph{strong derived length} of $P$.  It is clear that strong solvability implies solvability, but the converse is in general not true.
	
We say that a Poisson ideal $I$ of $P$ is \emph{nil} if all elements of $I$ are nilpotent. In particular, $I$ is said to be \emph{nil of bounded index} if there exists a positive integer $r$ such that $x^r=0$ for every $x\in I$. Furthermore, $I$ is said to be \emph{ associative nilpotent} if $I$ is a nilpotent ideal of $P$ as an associative algebra. 

\section{Lie properties of Poisson algebras}
 We start with the following result about strong solvability and strong Lie nilpotence, which represents the Poisson analogue of the corresponding results of Jennings for rings \cite[Theorems 5.5 and 6.5]{J}.
\begin{prop}\label{stronglynilp} Let $P$ be a  Poisson algebra over a field $\F$. Then the following statements hold.
\begin{enumerate}
\item If $P$ is strongly solvable, then $\ud_1(P)$  is associative
nilpotent of index at most $2^n-1$, where $n$ is the strong derived length of $P$.
\item $P$ is strongly Lie nilpotent if and only if $P$ is Lie nilpotent and strongly solvable.
\end{enumerate} 
\end{prop}
\Proof 
(1) It will be enough to show that, for every positive integer $k$, one has
\begin{equation}\label{sLieder}
\tilde{\delta}_1(P)^{2^k-1}\subseteq \tilde{\delta}_k(P).
\end{equation}
We proceed by induction on $k$. Let $k>1$, the claim being trivial for $k=1$. For all $x_1,x_2\in \ud_1(P)^{2^{k-1}-1}$ and $y_1,y_2\in P$ we have
 $$
x_1 x_2\{y_1,y_2\}=\{x_1 y_1,x_2 y_2\}-\{x_1,x_2 y_2\}y_1-\{x_1 y_1,x_2\}y_2+\{x_1,x_2\}y_1 y_2.
$$
It follows that
\begin{align*}
 \ud_1(P)^{2^{k}-1}&=\ud_1(P)^{2^{k-1}-1} \cdot \ud_1(P)^{2^{k-1}-1} \cdot \ud_1(P)\\
&\subseteq \{\ud_{k-1}(P), \ud_{k-1}(P)\}P= \tilde{\delta}_k(P),
\end{align*}
as claimed.

(2) Necessity is clear. Conversely, suppose that $P$ is both Lie nilpotent and strongly solvable. We proceed by induction on  the Lie nilpotence class $c$ of $P$. The claim is trivial for $c=1$.  Suppose then $c>1$. Then the Poisson algebra
$\bar{P}=P/\gamma_c(P)P$ is Lie nilpotent of class at most $c-1$ and so,  by the induction hypothesis, 
$\bar{P}$ is strongly nilpotent. Hence we have $P^{(r)}\subseteq \gamma_c(P)P$    for some positive integer $r$. As $\gamma_c(P) \subseteq Z_P(P)$, it follows that 
%
$$P^{(2r-1)}\subseteq \gamma_c(P)P^{(r)}\subseteq \gamma_c^2(P)P.$$ 
We can then continue in this fashion to show that 
\begin{equation}\label{gammac}
P^{(k(r-1)+1)}\subseteq \gamma_c^{k}(P)P
\end{equation}
for every nonnegative integer $k$.
 Now, by the first part of the theorem we know that $\ud_1(P)=\gamma_2(P)P$ is associative nilpotent. As $\gamma_c(P)\subseteq \ud_1(P)$, by Equation (\ref{gammac}) we conclude that $P^{(2k(r-1)+1)}=0$ for some $k$ sufficiently large, which proves the assertion. \qed

We now deal with solvability. Our aim is prove the Poisson version of a result obtained by Sharma and Srivastava \cite{SS} and, independently, by Smirnov and Zalesski \cite{SZ}.  
We will make use of the following result. Its proof can be found in \cite[Lemma 8.7]{MP} for $I=P$, by the same arguments actually work for any Lie ideal $I$ of the Poisson algebra $P$. 

\begin{lem}\label{MP} Let $P$ be a Poisson algebra and $I$ a Lie ideal of $P$. Then 
$$
\gamma_m(I)\gamma_n(I)\subseteq \gamma_{m+n-2}(I)P
$$
for every $m,n\geq 2$.
\end{lem}

Our first main result is the following:

\begin{thm}\label{S-Z} Let $P$ a solvable Poisson algebra of derived length $n$. Then the 
Poisson ideal generated by all elements $\{\{x_1,x_2\},\{x_3,x_4\},x_5\}$ is associative
nilpotent of index bounded by a function of $n$.
\end{thm}

\Proof
We can suppose $n\geq 2$, the claim being trivial when $P$ is abelian.  We have divided the proof into a sequence of steps.  Throughout the proof, $I$ is a Lie ideal of $P$.

{\bf Step 1}: $\gamma_3(I)\gamma_4(I)\subseteq \delta_2(I)P$. Let $x_1,x_2,x_3,x_4,x_5,x_6,x_7 \in I$.  We clearly have

$$
\{\{x_1,x_2\},\{x_3\{x_4,x_5\},x_6,x_7\} \}\in \{\gamma_2(I),\gamma_3(I)\}\subseteq \delta_2(L).
$$
On the other hand one has 
\begin{align*}
&\{\{x_1,x_2\},\{x_3\{x_4,x_5\},x_6,x_7\} \}=\{\{x_1,x_2\},\{x_3,x_6,x_7\}\{x_4,x_5\}+\\
&\{x_3,x_6\}\{x_4,x_5,x_7\}+\{x_3,x_7\}\{x_4,x_5,x_6\}+x_3\{x_4,x_5,x_6,x_7\} \}\\
 &\equiv \{\{x_1,x_2,x_3\} \cdot \{x_4,x_5,x_6,x_7\} \}\quad (\mathrm{mod}\, \delta_2(I)P),
\end{align*}
and the claim follows at once.

{\bf Step 2}: $\gamma_2(I)\{\gamma_2(I),P\}\subseteq \gamma_3(I)P$. Let $x_1,x_2,x_3,x_4\in I$ and $x_5\in P$. We have
\begin{align*}
\{x_3,x_4, \{x_1,x_2 x_5\}\}&=\{x_3,x_4, \{x_1,x_2\}x_5+x_2\{x_1,x_5\}\}\\
&\equiv \{ x_1,x_2\}\{x_3,x_4,x_5 \}\quad (\mathrm{mod}\, \gamma_3(I)P).
\end{align*} 
Since $\{x_3,x_4, \{x_1,x_2 x_5\}\}\in \gamma_3(I)P$, this yields the claim.

{\bf Step 3}: $\{I^2,P\}\subseteq I$. Let $x_1,x_2\in I$ and $x_3\in P$. We have
$$
\{x_1x_2,x_3\}=\{x_2,x_1x_3\}-\{x_2x_3,x_1\}\in I,
$$
yielding the claim.

{\bf Step 4}: $\gamma_2^2(I)\gamma_3(P)\subseteq I+\{\gamma_2(I),P\}P$. Let $x_1,x_2,x_3,x_4\in I$ and $x_5,x_6,x_7\in P$. Consider the element $x=\{x_1,x_2x_6,\{x_3x_5,x_4\},x_7\}\in \{\gamma_2(I),P\}$.  Then
$$
x=\{\{x_1,x_2\}x_6+x_2\{x_1,x_6\},\{x_3,x_4\}x_5+x_3\{x_5,x_4\},x_7\}.
$$
Since $x_2\{x_1,x_6\}$ and $x_3\{x_5,x_4\}$ are in $I^2$,  by Step 3 we have 
$$
\{\{x_1,x_2\}x_6,\{x_3,x_4\}x_5,x_7\}\in I+\{\gamma_2(I),P\}P.
$$
Moreover
$$
\{\{x_1,x_2, \{x_3,x_4\}x_5\}x_6,x_7\} \in \{\gamma_2(I),P\}P
$$
and so, by the Jacobi identity, we infer that
$$
\{\{x_1,x_2\}\{x_6,\{x_3,x_4\}x_5\},x_7\} \in I+\{\gamma_2(I),P\}P. 
$$
In turn, this implies
$$
\{\{x_1,x_2\}\{x_3,x_4\}\{x_6,x_5\},x_7\} \in I+\{\gamma_2(I),P\}P.
$$
It follows that
$$
\{x_1,x_2\}\{x_3,x_4\}\{x_5,x_6,x_7\} \in I+\{\gamma_2(I),P\}P,
$$
and the desired conclusion follows at once.

{\bf Step 5}: $\big(\gamma_2(I)^{5}\{\gamma_3(P),\gamma_3(P),P\} \big)^3\subseteq \delta_2(I)P$.
From Step 4 it follows that 
$$
\{\gamma_2(I)^2\gamma_3(P), \gamma_2(I)^2\gamma_3(P)\}\subseteq \gamma_2(I)+\{\gamma_2(I),P\}P.
$$
This forces
$$
\gamma_2(I)^4\{\gamma_3(P), \gamma_3(P)\}\subseteq \gamma_2(I)+\{\gamma_2(I),P\}P
$$
and then
\begin{equation}\label{ref6}
\gamma_2(I)^4\{\gamma_3(P), \gamma_3(P),P\}\subseteq \{\gamma_2(I),P\}P.
\end{equation}
Now, by Step 1 and Lemma \ref{MP} we obtain $\gamma_3(I)^3\subseteq \delta_2(I)P$ and so, by Step 2, we have  $\gamma_2(I)^3\{\gamma_2(I),P\}^3\subseteq \delta_2(I)P$. Therefore, in view of relation (\ref{ref6}) we get
$$
\gamma_2(I)^3\big(\gamma_2(I)^4\{\gamma_3(P),\gamma_3(P),P\} \big)^3\subseteq \delta_2(I)P,
$$
which is equivalent to the claim.

{\bf Step 6}: $\{\gamma_3(P),\gamma_2(P),P\}\{\gamma_2(P),\gamma_3(P),\gamma_3(P)\}\subseteq \{P,\gamma_3(P),\gamma_3(P)\}P$. Let $J,K$ be Lie ideals of $P$. Suppose that $K\subseteq L$, where $L$ is the Lie ideal consisting of all elements $x\in P$ such that $\{P,x\}\subseteq J$. We will show that 
\begin{equation}\label{ref8}
\{L,J,P\}\{P,P,K,K\}\subseteq \{J,K,P\}P.
\end{equation} 
From this, the desired conclusion will follow by setting  $J=K=\gamma_3(P)$ and $L=\gamma_2(P)$. 
Let $x_1\in L$, $x_2\in J$, $x_3,x_4,x_5\in P$, $x_6,x_7\in K$. We have 
\begin{align*}
&\{x_4x_6,x_7,x_2,x_3\}=\{\{x_4,x_7\}x_6+x_4\{x_6,x_7\},x_2,x_3\}\\
&=\{\{x_4,x_7,x_2\}x_6+\{x_4,x_7\}\{x_6,x_2\}+\{x_4,x_2\}\{x_6,x_7\}+x_4\{x_6,x_7,x_2\},x_3\}\\
&\equiv \{\{x_4,x_7,x_2\}\{x_6,x_3\}+\{x_4,x_7,x_3\}\{x_6,x_2\}+\{x_4,x_2,x_3\}\{x_6,x_7\}\\
&+\{x_4,x_2\}\{x_6,x_7,x_3\}+\{x_4,x_3\}\{x_6,x_7,x_2\}
\quad (\mathrm{mod}\,\{J,K,P\}P ).
\end{align*}
As $\{x_4x_6,x_7,x_2,x_3\}\in \{J,K,P\}$,  replacing $x_6$ with $\{x_5,x_6\}$ in the previous relation yields
\begin{align}\label{ref8a} \nonumber
\{x_4,x_7,x_2\}\{x_5,x_6,x_3\}&+\{x_4,x_7,x_3\}\{x_5,x_6,x_2\}\\
&+\{x_4,x_2,x_3\}\{x_5,x_6,x_7\}\in \{J,K,P\}P.
\end{align}
Now, by switching $x_6$ with $x_7$ and $x_4$ with $x_5$ in (\ref{ref8a}) we get
\begin{align*}
\{x_4,x_7,x_2\}\{x_5,x_6,x_3\}&+\{x_4,x_7,x_3\}\{x_5,x_6,x_2\}\\
&+\{x_5,x_2,x_3\}\{x_4,x_7,x_6\}\in \{J,K,P\}P
\end{align*}
and then
$$
\{x_4,x_2,x_3\}\{x_5,x_6,x_7\}\equiv \{x_5,x_2,x_3\}\{x_4,x_7,x_6\} \quad (\mathrm{mod}\,\{J,K,P\}P ),
$$
or equivalently
\begin{equation}\label{ref8b}
\{x_4,x_2,x_3\}\{x_6,x_5,x_7\}\equiv \{x_5,x_2,x_3\}\{x_7,x_4,x_6\} \quad (\mathrm{mod}\,\{J,K,P\}P ).
\end{equation}
Note that
\begin{equation}\label{8c}
\{x_5,x_6,x_4,x_7\}\equiv \{x_5,x_6,\{x_4,x_7\}\} \quad (\mathrm{mod}\,\{J,K,P\}P ).
\end{equation}
Therefore,  replacing $x_4$ by $x_1$ and $x_6$ by $\{x_6,x_4\}$  in (\ref{ref8b}) and using  (\ref{8c}) we get
\begin{align*}
&\{x_1,x_2,x_3\}\{x_6,x_4,x_5,x_7\}\equiv \{x_5,x_2,x_3\}\{x_7,x_1,\{x_6,x_4\}\}\\
&\equiv  \{x_5,x_2,x_3\}\{x_6,x_4,x_1,x_7\}\equiv  \{x_5,x_2,x_3\}\{x_6,x_1,x_4,x_7\}\\
&\equiv  \{x_5,x_2,x_3\}\{x_6,x_1,\{x_4,x_7\}\}\equiv  \{x_1,x_2,x_3\}\{x_4,x_7,x_5,x_6\}\\
&\equiv  \{x_1,x_2,x_3\}\{x_6,x_5,\{x_4,x_7\}\} \equiv  \{x_1,x_2,x_3\}\{x_6,x_5,x_4,x_7\} \quad (\mathrm{mod}\,\{J,K,P\}P ).
\end{align*}
Consequently, by the Jacobi identity we must have 
$$
\{x_1,x_2,x_3\}\{x_6,\{x_4,x_5\},x_7\}\in \{J,K,P\}P
$$ 
and then
$$\{x_1,x_2,x_3\}\{x_4,x_5,x_6,x_7\}\in \{J,K,P\}P,$$
from which (\ref{ref8}) follows at once.

{\bf Step 7}: \emph{Let $A,B,C$ be Lie ideals of $P$ such that $\{P,B\}\subseteq A$. Then $\{P,C\}\{A,B\}\subseteq \{A,C\}P$}. Let $x_1\in P$, $x_2\in A$, $x_3\in B$, $x_4\in C$. We have
$$
\{\{x_1x_2,x_3\},x_4\} \equiv \{x_1,x_4\}\{x_2,x_3\} \quad (\mathrm{mod}\, \{A,C\}P).
$$
Since $\{\{x_1x_2,x_3\},x_4\}\subseteq \{A,C\}$, the claim follows.

{\bf Step 8}: $\{\delta_2(P),P\}^3\gamma_3(P)^4\subseteq \{\gamma_3(P),\gamma_3(P),P\}P$.  By applying Step 7 for $A=\gamma_3(P)$, $B=\gamma_2(P)$ and $C=\{\gamma_3(P),\gamma_2(P)\}$ we get

$$
\{\gamma_3(P),\gamma_2(P), P\}\{\gamma_3(P),\gamma_2(P)\}\subseteq \{\gamma_2(P),\gamma_3(P), \gamma_3(P)\}P.
$$
Thus, by Step 6 we obtain
\begin{equation}\label{ref11}
\{\gamma_3(P),\gamma_2(P), P\}^2\{\gamma_3(P),\gamma_2(P)\}\subseteq \{\gamma_3(P),\gamma_3(P),P\}P.
\end{equation}
We have 
 \begin{align}\label{pezzo}\nonumber 
\{\delta_2(P),P\}&=\{\gamma_2(P),\gamma_2(P),P\}\\
&\subseteq \{\gamma_2(P),P,\gamma_2(P)\}=\{\gamma_3(P),\gamma_2(P)\}.
\end{align}
Moreover
\begin{align}\label{ref11a}\nonumber
\{\gamma_2(P),\gamma_2(P),\gamma_2(P)\}&\subseteq \{\gamma_2(P),\gamma_2(P),P,P\}\\ 
&\subseteq \{\gamma_2(P),P,\gamma_2(P),P\}=\{\gamma_3(P),\gamma_2(P), P\}.
\end{align}
Let $x_1,x_2,x_3,x_4 \in \gamma_2(P)$ and $x_5\in P$. We have
$$
\{x_1,x_2,\{x_3,x_4x_5\}\}\equiv \{x_1,x_2,x_5\}\{x_3,x_4\} \quad (\mathrm{mod}\, \{\gamma_2(P),\gamma_2(P),\gamma_2(P)\}P).
$$
As $\{x_1,x_2,\{x_3,x_4x_5\}\}\in \{\gamma_2(P),\gamma_2(P),\gamma_2(P)\}$, we have
$$ 
\{\gamma_2(P),\gamma_2(P),P\}\{\gamma_2(P),\gamma_2(P)\}\subseteq \{\gamma_2(P),\gamma_2(P),\gamma_2(P)\}P
$$
and so, by (\ref{ref11a}), we conclude that
\begin{equation}\label{ref11b}
\{\delta_2(P),P\}\delta_2(P)P\subseteq \{\gamma_3(P),\gamma_2(P), P\}P.
\end{equation}
Furthermore, for every $x_1,x_2\in \gamma_2(P)$ and $x_3,x_4\in P$ one has
$$
\{x_1,\{x_2x_3,x_4\}\}\equiv \{x_1,x_3\}\{x_2,x_4\}\quad (\mathrm{mod}\, \delta_2(P)P).
$$
Since 
$\{x_1,\{x_2x_3,x_4\}\}\in \delta_2(P)$, it follows that
\begin{align}\label{ref11c}
\gamma_3(P)^2\subseteq \delta_2(P)P.
\end{align}
At this stage, the combination of  (\ref{ref11}) with (\ref{pezzo}),  (\ref{ref11b}) and  (\ref{ref11c}) yields 
 $$
\{\delta_2(P),P\}^3\gamma_3(P)^4\subseteq \{\gamma_3(P),\gamma_3(P),P\}P.
$$
{\bf Step 9}:  \emph{We have $\{\delta_2(P),P\}^t=0$, where
$t=15^{n-2}+\frac{3}{2}(15^{n-2}-1).$} By Step 5 and Step 8 we infer that
$$
\gamma_2(I)^{15} \gamma_3(P)^{12}\{\delta_2(P),P\}^9\subseteq \delta_2(I)P.
$$
In particular, for $I=\delta_{i}(P)$ where $i\geq 0$ the previous relation becomes
\begin{equation}\label{passogen}
\delta_{i+1}(P)^{15} \gamma_3(P)^{12}\{\delta_2(P),P\}^9\subseteq \delta_{i+2}(P)P.
\end{equation}
Now, it follows from (\ref{passogen})  that
\begin{equation}\label{primopasso}
\delta_{n-1}(P)^{15} \gamma_3(P)^{12}\{\delta_2(P),P\}^9=0.
\end{equation}
Furthermore, by combining  (\ref{passogen}) (applied for $i=n-3$) and (\ref{primopasso}) we obtain
 $$
\delta_{n-2}(P)^{15^2} \gamma_3(P)^{12(15+1)}\{\delta_2(P),P\}^{9(15+1)}=0.
$$
At this stage, for every $j\geq n$ we see  by induction that 
$$
\delta_{n-j}(P)^{15^j} \gamma_3(P)^{12\cdot (\sum_{k=0}^{j-1}15^k)}\{\delta_2(P),P\}^{9\cdot (\sum_{k=0}^{j-1}15^k)}=0,
$$
thus
 \begin{equation}\label{passoj}
\delta_{n-j}(P)^{15^j} \gamma_3(P)^{\frac{6}{7}\cdot (15^j-1)}\{\delta_2(P),P\}^{\frac{9}{14}\cdot (15^j-1)}=0.
\end{equation}
Since $\{\delta_2(P),P\}\subseteq \delta_2(P) \subseteq \gamma_3(P)$, for $j=n-2$ the formula (\ref{passoj}) implies
 $$
\{\delta_2(P),P\}^{15^{n-2}+\frac{3}{2}(15^{n-2}-1)}=0,
$$
as claimed. This completes the proof. \qed



Let $P$ be a Poisson algebra over a field of characteristic zero. If $P$ is solvable, then it follows from \cite[Theorem 7.2]{MPR} that $\ud_1(P)$ is nil. Our next goal is to extend this result in  characteristic $p>2$. In this case, we will obtain the stronger conclusion that  $\ud_1(P)$ is nil of bounded index. Note that this is not true in characteristic zero.   Consider for instance the Poisson-Grassmann algebra ${\bf G}$ on a countable-dimensional vector space $V=\mathrm{span} \{e_1,e_2,\ldots\}$ over a field  of characteristic zero (cf. \cite[\S 6]{MPR}). Then  ${\bf G}$ is Lie nilpotent of class 2 (see \cite[Theorem 6.1(5)]{MPR}). However, 
$$
\{{\bf G},{\bf G}\}=\{e_{i_1}e_{i_2}\cdots e_{i_{2n}}\vert\, i_1<i_2<\cdots<i_{2n}\}
$$
is not nil of bounded index since, for every positive integer $n$, the element $e_1e_2+e_3e_4+\cdots+e_{2n-1}e_{2n}$ has index of nilpotency $n+1$. 

We first prove a couple of preliminary lemmas.

\begin{lem}\label{x,y}
	 Let $I$ be a Poisson ideal of a   Poisson algebra  $P$ over a field $\F$.  Then, for every  $a,b,c,x,y\in I$ and
	 positive integer $m$, we have 
	 \begin{enumerate}
	 	\item $	(\{a, b, c\}-\{b, c, a\})\gamma_m(I)\su \gamma_{m+1}(I)P$.
	\item $\{x,y\}^2\gamma_m(I)\su \gamma_{m+1}(I)P$.
	 \end{enumerate} 
\end{lem}
\Proof Let $d, e\in I$. Direct computations imply that
\begin{align*}
(\{a, b, c\}-\{b, c, a\})\{d, e\}=& \{ac, d, b, e\} -a\{c, d, b, e\} +\{c, d, be, a\}\\
& +b\{d, c, e, a\} +\{d, c, b,  a\}e + \{d,a,b,c\}e \\
&-\{a,d,b,e\}c+\{b,c , ae, d\}+ a\{c,b,e,d\}\\
&+\{c,b,a,d\}e+\{b,a,ec,d\}+e\{a,b,c,d\}\\
&+\{a,d,be,c\}+\{a,b,e,d\}c-b\{a, d, e, c\}.
\end{align*}
It is now enough to take  $d\in \gamma_{m-1}(I)$ to deduce Part (1).

By Part (1), for all $a,b\in I$, we have 
\begin{equation}\label{xyy}
\{a,b,b\}\gamma_m(I)\su \gamma_{m+1}(I)P.
\end{equation}
Using Part (1) again, we get  $\{a, b, c\}z\equiv \{b, c, a\}z\; (\text{mod } \gamma_{m+1}(I)P)$, for every $a, b, c\in I$ and  $z\in \gamma_m(I)$.
It follows from the Jacobi identity that 
\begin{align}\label{abc}
3\{a, b, c\}\gamma_m(I)\su \gamma_{m+1}(I)P.
\end{align}
Let $x,y,z\in I$. A direct expansion yields
\begin{align}\label{ab}
\{x,y\}\{y,z\}&=3\{xy,y,z\}-3\{x,y,z\}y\nonumber\\
&-\{xz,y,y\}+\{x,y,y\}z+x\{z,y,y\}.
\end{align}
The result now follows  from (\ref{xyy}),  \eqref{abc}, and \eqref{ab}.  \qed

We note that Equation \eqref{abc} immediately implies that if $P$ is Lie nilpotent and $\text{char}(\F)\neq 3$, then $\gamma_{3}(P)P$ is associative nilpotent. This also follows from  \cite[ Lemma 8.7]{MP}.

We also need the following  technical lemma after which we would be ready to prove the second main result of this paper.

\begin{lem}\label{delta2J}
	Let $P$ be a Poisson algebra and denote by $\mathcal{J}$  the Poisson ideal of $P$ generated by all elements $\{\{x_1, x_2\}, \{x_3, x_4\}, x_5\}$ with $x_1,\ldots ,x_5 \in P$. Then $4\{\{x_1, x_2\}, \{x_3, x_4\}\}^3\in   \mathcal{J}$,  for all $x_1,\ldots, x_4\in P$.
\end{lem}
\Proof Let $z_1=\{x_1, x_2\}$ and $z_2=\{x_3, x_4\}$. Note that 
\begin{align*}
2\{z_1, z_2\}^2=& \{z_1^2, z_2, z_2\}- 2z_1\{z_1, z_2, z_2\}\equiv \{z_1^2, z_2, z_2\}\quad (\mathrm{mod}\, \mathcal{J}).
\end{align*}
Furthermore, 
\begin{align*}
2\{z_1^2, z_2, z_2\} \{z_1, z_2\}=& \{z_2^2, z_1^2, z_2, z_1\}+2z_2\{z_1^2, z_2, z_2, z_1\} \equiv 0 \quad (\mathrm{mod}\, \mathcal{J}).
\end{align*}
We conclude that $4\{z_1,z_2\}^3\in   \mathcal{J}$, as required. \qed

\begin{thm}\label{solvablenil}  Let $P$ be   Poisson algebra  over a field $\F$ of characteristic $p\geq 0$. 
	\begin{enumerate}
		\item If $P$ is solvable and $p\neq 2$, then $\{P, P\}P$ is a  nil ideal.
		\item  If $P$ is solvable and $p\geq 3$, then $\{P, P\}P$ is  nil of bounded index.
		\item   If $P$ is Lie nilpotent and $p>0$, then $\{P, P\}P$  is  nil of bounded index.
	\end{enumerate}  
\end{thm}
\Proof   Let $x,y\in P$. Note that, by Lemma \ref{x,y},  $\{x,y\}^3\in \gamma_{3}(P)P$ and  $\{x,y\}^5\in \gamma_{4}(P)P$. As a consequence, by Step 1 of the proof of Theorem \ref{S-Z}, we deduce that  $\{x,y\}^8\in \delta_{2}(P)P$. 
Now, a typical element $z$ of $\{P, P\}P$  is  of the form 
$z=\sum_{i=1}^s \{x_i, y_i\}z_i,$ where $x_i,y_i,z_i\in P$ for every $1\leq i\leq s$. 

To prove Part (1), we note that  
 $z^{8s}\in \delta_{2}(P)P$. Therefore it suffices to show that $\delta_{2}(P)P$ is nil.   
Now, every element $a$ of $\delta_2(P)P$ is  of the form 
$a=\sum_{i=1}^r\{u_i,v_i\}x_i$, where $u_i, v_i\in \{P, P\}$ and $x_i\in P$, for every $1\leq i\leq r$.
Since $p\neq 2$, it follows from Lemma \ref{delta2J} that  each $\{u_i,v_i\}^3\in \mathcal{J}$. Thus  $a^{3r}\in \mathcal{J}$.   Theorem \ref{S-Z}  now implies that $a$ is nilpotent. 

Suppose that $p\geq 3$ and let us prove Part (2). For what was showed above, for all $x, y\in P$ we have $\{x,y\}^8\in \delta_{2}(P)P$, so that $\{x,y\}^{p^2}\in \delta_{2}(P)P$. Hence $z^{p^2}=\sum_{i=1}^s \{x_i, y_i\}^{p^2}z_i^{p^2}\in \delta_{2}(P)P$. 
It is now enough to observe that  $\delta_{2}(P)P$ is nil of bounded index. Indeed, every element $b$ of $\delta_2(P)P$ is of the form 
$b=\sum_{i=1}^r\{u_i,v_i\}x_i$, where $u_i, v_i\in \{P, P\}$ and $x_i\in P$, for every $1\leq i\leq r$.
Since $p\neq 3$, it follows from Lemma \ref{delta2J} that  each $\{u_i,v_i\}^p\in \mathcal{J}$. Hence  $b^{p}\in \mathcal{J}$, and by  Theorem \ref{S-Z}   we obtain the desired conclusion. 

To prove Part (3), it is enough  by Part (2)  to prove the statement for  $p=2$. By Lemma \ref{x,y} and  the Jacobi identity, we have
 \begin{align}\label{abc2}
 \{a, b, c\}\gamma_m(P)\su \gamma_{m+1}(P)P,
 \end{align}
for every  $a, b, c\in P$. Let $s$ be the nilpotence class of $P$ and $r$ the smallest integer such that $2^r\geq s+1$.
 We deduce from Equation \eqref{abc2} that $z^{2^r}=0$, for every $z\in  \gamma_{3}(P)$. We can now replace $P$ with $P/\gamma_{3}(P)P$ and assume that $\gamma_{3}(P)=0$. It then follows from Lemma \ref{x,y} that 
 $\{x,y\}^4=0$, for every $x, y\in P$. Since every element $z\in \{P, P\}$ is a linear combination of 
 the $\{x_i,y_i\}$'s, it then follows that $z^4=0$. This finishes the proof.
\qed

\begin{rem} \emph{The assumption on the characteristic of the ground field cannot be dropped in Theorem \ref{solvablenil}(1). In fact, the Hamiltonian Poisson algebra ${\bf H}_2$ in two indeterminates over a field $\F$ of characteristic 2 is solvable of derived length 3 (see \cite[Lemma 11.3]{MP}). However, as ${\bf H}_2$ is a domain, $\tilde{\delta}_1({\bf H}_2)$ is not nil.
}
\end{rem}

Let us now discuss some consequences of our previous results. 
We say that a Poisson algebra $P$ is \emph{reduced} if $P$ is reduced as an associative algebra (that is, $P$ is free of nonzero nilpotent elements). By Theorem \ref{solvablenil} and Proposition \ref{stronglynilp} we have  

\begin{cor}\label{reduced} Let $P$ be a reduced Poisson algebra over a field $\F$. Then the following statements hold. 
\begin{enumerate}
\item $P$ is Lie nilpotent if and only if $P$ is abelian.
\item $P$ is strongly solvable if and only if $P$ is abelian.
\item If $\F$ has characteristic not 2, then $P$ is solvable if and only if $P$ is abelian.
\end{enumerate}
\end{cor}

Let $L$ be a Lie algebra $L$ over a field $\F$ and denote by $S(L)$ its symmetric algebra, which we identify with the polynomial ring $\F[x_1,x_2,\ldots]$ where $x_1,x_2,\ldots$ is an $\F$-basis of $L$ over $\F$. The Lie bracket $\{x,y\}$ of $L$ can be uniquely extended to a Poisson bracket of  $S(L)$ so that this commutative algebra becomes a Poisson algebra, called the \emph{symmetric Poisson algebra} of $L$. 
The Lie identities of $S(L)$ have been investigated by Monteiro Alves and Petrogradsky in \cite{MP} (see also \cite{S} for further developments). In particular, by generalizing a result of Shestakov in \cite{Sh}, in \cite[Theorem 4.4]{MP} the authors proved the theorem quoted below:

\begin{thm}
Let $L$ be a Lie algebra over a field $\F$, and $S(L)$ its symmetric Poisson algebra. The following conditions are equivalent: 
\begin{enumerate}
\item $L$ is abelian;
\item $S(L)$ is strongly Lie nilpotent;
\item $S(L)$ is Lie nilpotent;
\item  $S(L)$ is strongly solvable;
\item   $S(L)$  is solvable (here assume that $\charac (\F) \neq 2$).
\end{enumerate}
\end{thm}
 As the Poisson algebra $S(L)$ is reduced, the previous theorem is now a consequence of Corollary \ref{reduced}.

When $\F$ has characteristic $p>0$, the Poisson bracket of $S(L)$ naturally induces a Poisson bracket on the factor algebra ${\bf s}(L)=S(L)/I$, where $I$ is the ideal generated by the elements $x^p$ with $x\in L$. The conditions under which the Poisson algebra ${\bf s}(L)$ is solvable in characteristic different from 2 or Lie nilpotent were also determined by  Monteiro Alves and Petrogradski in \cite{MP}. In characteristic 2, solvability of $S(L)$ and ${\bf s}(L)$ was left by the authors as an open problem and a related conjecture was proposed (see \cite[\S 5.3, Conjecture]{MP}). A counterexample for this conjecture was found in \cite{S}.  However, by means of Theorem \ref{S-Z}, a corrected version of  \cite[\S 5.3]{MP} can be proved, thereby completing the classification. This will be accomplished in a  forthcoming separate paper.

\end{document}